\documentclass[10pt]{article}
\usepackage[dvips]{graphicx}
\usepackage{amssymb}
\usepackage{amsmath}

\usepackage{caption}
\usepackage{hyperref}
\usepackage{arcs}
\usepackage {xcolor} 
\usepackage{subfigure}
\usepackage{authblk}
\usepackage{calrsfs}

\newtheorem{theorem}{Theorem}
\newtheorem{Proposition}{Proposition}

\newtheorem{Remark}{Remark}

\usepackage{graphicx} 

\newenvironment{AMS}{\small\bf 2020 AMS subject classification: }{} 

\title{A note on polynomial-free unisolvence of polyharmonic splines at random points}

\author[1]{Len Bos}
\affil[1]{University of Verona, Italy}
\author[2]{Alvise Sommariva} 
\affil[2,3]{University of Padova, Italy}
\author[3]{Marco Vianello}

\date{\today}

\begin{document}

\maketitle

\begin{abstract}
In this note we prove almost sure unisolvence of RBF interpolation on randomly distributed sequences by a wide class of polyharmonic splines (including Thin-Plate Splines), without polynomial addition.
\end{abstract}

\vskip0.2cm
\noindent
\begin{AMS}
{\rm 65D05,65D12.}
\end{AMS}
\vskip0.2cm
\noindent
{\small{\bf Keywords:} multivariate interpolation, Radial Basis Functions, polyharmonic splines, unisolvence.}

\vskip1cm
Interpolation by RBF (Radial Basis Functions) is nowadays one of the basic tools of computational mathematics, with a special relevance in the framework of multivariate approximation by scattered data and of meshless methods for numerical modelling. In the literature, since the beginning of RBF theory, unisolvence has been proved mainly resorting to positive definiteness, which becomes {\em conditional} for some RBF such as multiquadrics and polyharmonic splines (CPD, Conditionally Positive Definite of order $m$); cf. e.g. \cite{F07,I03,W05} for definitions and properties. 

On the other hand, while for multiquadrics and other CPD of order $m=1$ it is known that the interpolation matrix at distinct nodes is nonsingular \cite{M86}, unisolvence without the addition of a polynomial term of degree $m-1$ seems not investigated for other popular CPD RBF, such as polyharmonic splines. 
We recall that polyharmonic  splines correspond to the radial functions $$\phi(r)=r^{2k} \log(r)\;,\;\;k\in \mathbb{N}$$ (TPS, Thin-Plate Splines, order $m=k+1$) and $$\phi(r)=r^\nu\;,\;\;0<\nu\notin 2\mathbb{N}$$ (RP, Radial Powers, order $m=\lceil\nu/2\rceil$);
cf. e.g. \cite{D76,F07,I04,Meing84}. One of their most relevant features is that they are scale invariant, namely the condition number of the interpolation matrix $[\phi(\varepsilon\|x_i-x_j\|_2)]$ and the cardinal basis functions of interpolation
are independent of the scale parameter $\varepsilon>0$; cf. e.g. \cite{I03}. We recall that with scale dependent bases the scaling choice and effect are a delicate matter and still an active research topic in RBF interpolation; cf. e.g. \cite{LS23} with the references therein. 

The following statement appeared in the treatise \cite{F07}: {\em ``There is no result that states that interpolation with Thin-Plate Splines (or any other strictly conditionally positive definite function of order $m\geq 2$) without the addition of an appropriate degree $m-1$ polynomial is well-posed''}. This is not completely true, since there was a known result for univariate RP with $\nu=3$ ($m=2$), cf. \cite{BS87}. The only other known result concerns the case $m=1$, cf. \cite{M86}.
To our knowledge, the situation does not seem to have changed thenceforth. 

In this note we prove that interpolation at {\em random points} by TPS or RP with noninteger exponent is {\em almost surely unisolvent} without the need of adding a polynomial term. We shall use the key property of TPS and RP of being {\em analytic functions} up to their center, where they present a singularity of some derivative. 
Recently, the  nonstandard case of fixed centers (different from the sampling points) has been treated in \cite{DASV23} (we may also mention \cite{XN23} where however only spaces of everywhere analytic functions are considered).

Here we cope the more difficult case of classical RBF interpolation where the centers and the sampling points coincide. 
The fact that unisolvence in practice occurs for TPS without polynomial addition has been clearly recognized and partially discussed in the recent paper \cite{P22}, though only via numerical results.

\begin{theorem} 
Let $\Omega$ be an open connected subset of $\mathbb{R}^d$, $d\geq 2$, and 
$\{x_i\}_{i\geq 1}$ be a randomly distributed sequence on $\Omega$ with respect to any given probability density $\sigma(x)$, i.e. a point sequence produced by sampling a sequence of absolutely continuous random variables $\{X_i\}_{i\geq 1}$ which are independent and identically distributed in $\Omega$ with density $\sigma\in L^ 1_+(\Omega)$. Moreover, let $V_n=[\phi_j(x_i)]$, $\phi_j(x)=\phi(\|x-x_j\|_2)$, $1\leq i,j\leq n$, be the interpolation matrix with respect to TPS or RP with noninteger exponent. 

Then, for every $n\geq 1$ the matrix $V_n$ is {\em a.s.}\ (almost surely) nonsingular. 
\end{theorem}

\vskip0.3cm
\noindent{\bf Proof.} Before starting the proof, some observations are in order. Notice that the diagonal of $V_n$ is zero, since with polyharmonic splines $\phi(0)=0$. Moreover, 
in any dimension $d\geq 2$ the functions $\{\phi_j(x)\}$ are linearly independent in $\Omega$ iff the centers are distinct (which is a necessary condition for unisolvence). In fact, each $\phi_j$ is analytic in $\Omega\setminus \{x_j\}$ and singular at $x_j$, by analyticity of $\phi(\sqrt{\cdot})$ in $\mathbb{R}^+$, being composition of analytic functions singular at 0 like  $\log(\sqrt{\cdot})$ and $(\sqrt{\cdot})^\nu$, $0<\nu\notin 2\mathbb{N}$; cf. e.g. \cite{KP02}. If the $\{\phi_j(x)\}$ were linearly dependent with distinct centers, one of them being linear combination of the others would become analytic at its own center. Finally, we recall the basic property that a subset of $\Omega$ has null measure with respect to $d\sigma=\sigma(x)\,dx$, if it has null Lebesgue measure. A set with null Lebesgue measure is usually called a ``null set'' in measure theory. The condition becomes also necessary when $\sigma$ is almost everywhere positive.  

First, we prove that $x_1$ and $x_2$ are a.s. distinct and $V_2$ is a.s. nonsingular. Indeed, $det(V_2)=\phi_2(x_1)\phi_1(x_2)=\phi^2(\|x_1-x_2\|_2)=0$ iff $x_1=x_2$. But the probability that $x_2=x_1$ given $x_1$ is null, because in general the probability that a random point belongs to a finite set is null, a finite set being a null set. Next, we prove that the same properties holds for $n=3$. The fact that $x_1$, $x_2$ and $x_3$ are a.s. distinct follows from the same argument above. On the other hand, it is a general fact that $3\times 3$ matrices with null diagonal and positive extra-diagonal entries have a positive determinant, as it can be easily checked. Indeed, it turns out that for any such $A=[a_{ij}]$, $1\leq i,j\leq 3$, we have $det(A)=a_{12}a_{23}a_{31}+a_{13}a_{21}a_{32}>0$. In our case the property holds a.s., since the points are a.s. distinct.

Now, let us assume almost sure nonsingularity for $n\times n$ interpolation matrices on random sequences and that the points $x_1,\dots,x_n$ are a.s distinct with $n\geq 3$, and prove that the same properties hold with $n+1$. Consider the $(n+1) \times (n+1)$ matrix $U_{n+1}(x)$ obtained by adding to $V_n$ the $(n+1)$-th column $[\phi_1(x),\dots,\phi_n(x),0]^t$ and the $(n+1)$-th row 
$[\phi_1(x),\dots,\phi_n(x),0]\;,$
$$
U_{n+1}(x)=\left(\begin{array} {cc}
V_{n} & \vec{\Phi}_n^t(x)\\
\\
\vec{\Phi}_n(x) & 0
\end{array} \right)
\;,\;\;
\vec{\Phi}_n(x)=[\phi_1(x),\dots,\phi_{n}(x)]\;
$$

Applying Laplace rule to the $(n+1)$-th row, we get 
$$
f_n(x)=det(U_{n+1}(x))=\alpha_n(x)\phi_n(x)+\dots
+\alpha_1(x)\phi_1(x)\;,
$$
where $\alpha_1,\dots,\alpha_n$ are the corresponding minors with the appropriate sign. 
Moreover, developing the minor 
corresponding to $\alpha_n(x)$ by its last row 
and putting in evidence the factor $\phi_n(x)$, we get 
$$
f_n(x)=d_{n-1}\phi_n^ 2(x)+A_n(x)\phi_n(x)+B_n(x)\;,\;\;d_{n-1}=-det(V_{n-1})\;,
$$
where $A_n(x)\in span(\phi_1,\dots,\phi_{n-1})$ and $B_n(x)\in span(\phi_j\phi_k\,,\,1\leq j,k\leq n-1)$. 
Then $f_n$ is a.s. not identically zero in $\Omega$. 

Indeed, let us take a neighborhood of $x_n$ where all the $\phi_j$, $1\leq j \leq n-1$, are analytic, and fix all the variables but one to the corresponding coordinates of $x_n$. Then the problem is reduced to show that the univariate function $g(t)=\phi(|t-t_0|)=(t-t_0)^{2k}\log(|t-t_0|)$ cannot satisfy the identity $g^2\equiv\alpha g+\beta$, where $\alpha(t)$ and $\beta(t)$ are analytic in a neighborhood of $t_0$ (where $t_0$ is the coordinate of $x_n$ corresponding to the free univariate variable, say $t$).

A first proof in the case of TPS is the following. Observe that $\alpha\not\equiv 0$, otherwise $g^2\equiv \beta$ would be analytic, whereas it has a singularity of the $4k$-th derivative at $t_0$. In fact, taking for simplicity $t> t_0$ and applying Leibniz rule for the derivatives of a product to $g^2$, we get that the lower derivatives of $g^2$ are regular whereas $D^{4k}g^2$ has a $\log^2$ singularity at $t_0$ (and the same holds for $\alpha g\equiv g^2-\beta$). Since $\alpha\not\equiv 0$, the Taylor series of $\alpha$ centered at $t_0$ cannot have all zero coefficients, i.e. there exists a minimum $m\geq 0$ such that  $D^m \alpha(t_0)\neq 0$. 

Then, $\alpha(t)=\sum_{j=m}^\infty D^j \alpha(t_0)(t-t_0)^j/j!$, and for $t>t_0$
$$\alpha(t)g(t)=\sum_{j=m}^\infty D^j \alpha(t_0)(t-t_0)^{j+2k}\log(t-t_0)/j!
$$
so that, again by Leibniz rule applied termwise to the series, the lower derivatives of $\alpha g$ are regular whereas $D^{m+2k}(\alpha g)$ has a $\log$  singularity at $t_0$, which gives a contradiction for any value of $m$. 

The reasoning just developed does not seem to extend in the case of RP with noninteger exponent. 
We give now a second proof, which is valid for both, TPS and RP with noninteger exponent. Indeed, if $g^2\equiv\alpha g+\beta$ then $g(t)$ solves for every fixed $t$ the quadratic equation $g^2(t)-\alpha(t)g(t)-\beta(t)=0$, that is $g(t)$ necessarily takes one of the values $$g_\pm(t)=\frac{1}{2}\left(\alpha(t)\pm \sqrt{\Delta(t)}\right)\;,$$ where 
$\Delta(t)=\alpha^2(t)+4\beta(t)$ is the discriminant. 

First, $\Delta(t_0)=0$, otherwise the sign would be preserved in a neighborhood of $t_0$. If $\Delta(t_0)>0$ then by continuity $g\equiv g_+$ or alternatively $g\equiv g_-$ there, with both analytic at $t_0$, whereas $g$ is singular at $t_0$ having a singularity of $D^{2k}$ for TPS and of $D^{\lceil\nu\rceil}$ for RP. If on the contrary $\Delta(t_0)<0$, both $g_+(t_0)$ and $g_-(t_0)$, and hence $g(t_0)$ would be complex nonreal at $t_0$, whereas $g(t_0)=0$. 

Now, since $\Delta(t)$ is analytic, by a well-known result $\Delta(t)=(t-t_0)^mh(t)$, where $m$ is called the order of the zero, $1\leq m\in \mathbb{N}$, and $h$ is (locally) analytic and nonvanishing at $t_0$ (cf., e.g., \cite{KP02}). Then, if $h(t_0)<0$ we have $h(t)<0$ in a neighborhood and $\Delta(t)<0$ for $t>t_0$ (any $m$), i.e. $g_\pm(t)$ and hence $g(t)$ are complex nonreal for $t>t_0$. 

If $h(t_0)>0$ we have $h(t)>0$ in a neighborhood. For $m$ odd, $\Delta(t)<0$ for $t<t_0$, i.e. $g_\pm(t)$ and hence $g(t)$ are complex nonreal for $t<t_0$. For $m$ even, by continuity $g\equiv g_+$ or alternatively $g\equiv g_-$ for $t>t_0$, we get $g_\pm(t)=\frac{1}{2}\left(\alpha(t)\pm (t-t_0)^{m/2}\sqrt{h(t)}\right)$. Notice that all the derivatives of  $g_\pm(t)$ have a finite limit for $t\to t_0^+$. 
In all cases, we get a contradiction for TPS and RP with noninteger exponent, since $g(t)$ is real (positive) for $t\neq t_0$, and singular at $t_0$ with a suitable derivative going to $\infty$ for $t\to t_0$.

 We can now use the fact that $f_n(x)$ is analytic in the open connected set $\Omega\setminus \{x_1,\dots,x_n\}$ 
 (notice that the latter has a unique connected component in dimension $d\geq 2$), being sum of products of analytic functions, and is a.s. not identically zero there, otherwise by continuity it would be a.s. identically zero on the whole $\Omega$, which has been just excluded. 

Clearly, the points $\{x_1,\dots,x_n,x_{n+1}\}$ also are a.s. distinct because such are $\{x_1,\dots,x_n\}$ and the probability that $x_{n+1}$ coincides with one of them is null, a finite set being a null set. 
On the other hand, 
$V_{n+1}=U_{n+1}(x_{n+1})$ holds since $\phi_{n+1}(x_{n+1})=0$ and $\phi_j(x_{n+1})=\phi_{n+1}(x_{j})$ for $j=1,\ldots,n$. 

Then, 
$det(V_{n+1})=det(U_{n+1}(x_{n+1}))=f_n(x_{n+1})$ is a.s. nonzero, 
since the zero set of $f_n$ is a null set by a well-known basic result of measure theory, asserting that the zero set of a not identically zero real analytic function on an open connected set in $\mathbb{R}^d$ is a null set (cf. \cite{M20} for an elementary proof). \hspace{0.2cm} $\square$

\begin{Remark}
{\em We recall that that for the RP with $\nu=1$ unisolvence is well-known, due to invertibility of distance matrices for any choice of distinct points in any dimension (cf. \cite{M86}). The same can be said for $\nu=3$ in dimension $d=1$, by the result proved in \cite{BS87}. To our knowledge, nothing is known for integer $\nu>3$, for non integer $\nu>0$ and for TPS. 
On the other hand, 
the mere fact that the points are distinct does not guarantee unisolvence for TPS. For example, if 
$n-1$ points are taken, deterministically or even randomly, on a unit sphere (the boundary of a unit ball) centered at another fixed one, 
the resulting $n\times n$ TPS interpolation matrix has a null row (since $\log(1)=0)$ and hence is singular.

We observe finally that the present proving approaches do not seem to work for univariate instances (where $\Omega\setminus \{x_1,\dots,x_n\}$ is disconnected), as well for multivariate RP with odd integer exponent. These cases might deserve further investigations.
}
\end{Remark}

\begin{Remark}
{\em We stress that our result does not mean that enriching the functional space by polynomials, as in classical CPD RBF interpolation, is useless, because requiring exactness on low-dimensional polynomial spaces could be driven by specific application purposes. For example the need of recovering exactly constant and linear polynomials can arise in meshless methods for PDEs, cf. e.g. \cite{F07}. Moreover, working with positive-definite symmetric matrices could have computational advantages with respect to simply nonsingular symmetric matrices. 
Indeed, the present result only means that polynomial addition is {\em not necessary} within a wide class of polyharmonic splines to get unisolvence in the continuous random setting with respect to any density, a fact apparently unproved in the RBF literature till now. As noted in \cite{P22}, there are applications where polynomial addition is not really needed, for example point cloud modelling from high density geodetic data. In these cases, the result proved in the present note can be meaningful.}
\end{Remark}

\section*{Acknowledgements}

Work partially
supported by the
DOR funds 
of the University of Padova, and by the INdAM-GNCS.
This research has been accomplished within the RITA ``Research ITalian network on Approximation", the SIMAI Activity Group ANA\&A and the UMI Group TAA ``Approximation Theory and Applications".

\end{document}